\documentstyle[elsart12]{elsart}

\input amssym.def
\input amssym

\def\F{{\cal F}}
\def\Q{{\cal Q}}

\begin{document}
\begin{frontmatter}
\title{{\Large \bf Generic Gaussian ideals}}
\author[Rutgers]{Alberto Corso\thanksref{CNR}},
\thanks[CNR]{Corresponding author. Present address: Department of
Mathematics, Purdue University, West Lafayette, IN 47907, USA.
E-mail: corso@math.purdue.edu.}
\author[Rutgers]{Wolmer V. Vasconcelos\thanksref{NSF}},
\thanks[NSF]{E-mail: vasconce@rings.rutgers.edu.}
\author[Mexico]{Rafael H. Villarreal\thanksref{SNI}}
\thanks[SNI]{E-mail: vila@esfm.ipn.mx.}
\address[Rutgers]{Department of Mathematics,
Rutgers University,
New Brunswick, NJ 08903, USA}
\address[Mexico]{Instituto Polit\'ecnico Nacional,
Escuela Superior de F\'{\i}sica y Matem\'aticas,
07300 Mexico, D.F. MEXICO}
\date{December 5, 1995}

%
%\begin{keyword}
%Content of a polynomial, Gauss Lemma, Dedekind--Mertens
%formula, Cohen--Macaulay rings, Gorenstein ideals, Noether
%normalizations.
%\end{keyword}
%

\begin{abstract}
\noindent The content of a polynomial $f(t)$ is the ideal generated by
its coefficients. Our aim here is to consider a beautiful formula
of Dedekind--Mertens on the content of the product of two polynomials,
to explain some of its features from the point of view of Cohen--Macaulay
algebras and to apply it to obtain some Noether normalizations of certain
toric rings. Furthermore, the structure of  the primary decomposition of
generic products is given and some extensions to joins of toric rings
are considered.

\

\par\noindent {\it 1991 Math. Subj. Class.:} Primary 13H10; Secondary
13D40, 13D45, 13H15
\end{abstract}
\end{frontmatter}

\section{Introduction}
If $R$ is a commutative ring and  $f=f(t) \in R[t]$ is a polynomial,
say $f = a_0 + \cdots + a_mt^m$, the {\em content} of $f$ is the
$R$-ideal $(a_0,\ldots, a_m)$. It is denoted  by $c(f)$.
Given another polynomial $g$, the {\em Gaussian ideal} of $f$
and $g$ is the $R$-ideal
\begin{eqnarray}\label{Gaussideal}
G(f, g)= c(f g).
\end{eqnarray}

This ideal bears a close relationship to the ideal $c(f)c(g)$, one
aspect of which is expressed in the classical lemma of Gauss: If $R$ is
a PID then
\begin{eqnarray} \label{Gausslemma}
c(f g)= c(f)c(g).
\end{eqnarray}

In fact, if $R$ is a domain, then this equality holds for arbitrary
pairs of polynomials if and only if $R$ is a Pr\"ufer domain. In general,
these two ideals are very different but one aspect of their relationship
is given by
(see \cite{Northcott})
\begin{equation}\label{contentformula}
c(f g)c(g)^m = c(f) c(g)^{m+1}.
\end{equation}

One of our purposes in this note is to `explain' this formula,
originally due to Dedekind--Mertens, in terms of the theory of
Cohen--Macaulay rings, and to consider some extensions of it. More
precisely, we study the ideal $G(f,g)$ in the case when $f$ and $g$
are generic polynomials. It turns out that several aspects of the
theory of Cohen--Macaulay rings---e.g., $a$-invariants and linkage
theory---show up very naturally when we closely examine $G(f, g)$.

%
%The more usual formula (Lemma of Gauss),
%valid for domains with the $\gcd$ property, is the
%equality
%$C(f\cdot g) = C(f)\cdot C(g)$, where $C(f)$ is  $\gcd(c(f))$. It is
%obtained by cancelling a factor
%$(C(g))^d$ from both sides of the $C(\cdot)$--version
%of  (\ref{contentformula}).
%(Studies of cases when  $c(fg)=c(f)\cdot
%c(g)$ are made in \cite{GV} and \cite{HeHu}.)
%For more general rings this cancelling out is not always possible,
%and here we move in  the opposite direction.
%

One path to our analysis and its applications to Noether
normalizations of some semigroup rings starts by multiplying both
sides of (\ref{contentformula}) by $c(f)^m$; we obtain
\begin{equation}\label{contentformula2}
c(f g)[c(f) c(g)]^m = c(f) c(g) [c(f) c(g)]^{m}.
\end{equation}
It is this `decayed' content formula that will be the focus of our
observations. One result (namely, Theorem~\ref{sharpind}) will show that
(\ref{contentformula2}) is sharp in terms of the exponent $m
=\deg f$ (and, therefore, (\ref{contentformula}) as well).
It will be the outgrow of looking for Noether
normalizations of certain rings generated by monomials and basic
facts of the theory of Cohen--Macaulay rings. In particular
(\ref{contentformula2}) is shown to be a direct consequence of the
lemma of Gauss.

To make this connection, we recall the notion of a {\em reduction}
of an ideal (see \cite{NR}).
Let $R$ be a commutative Noetherian ring and let $I$ be an ideal. A
{\em reduction} of $I$ is an ideal $J \subset I$
such that, for some non-negative integer $r$, the equality $I^{r+1}=
J I^r$ holds. The smallest such integer is the {\em reduction
number} $r_J(I)$ of $I$ relative to $J$.
Thus (\ref{contentformula2}) says\footnote{See also \cite[Section
3]{RR}.} that  $J=c(f g)$ is a reduction for $I = c(f) c(g)$, and that
the reduction number is at most $\min\{\deg f, \deg g\}$.

One of the advantages of reductions is that they contain much of
the information carried by $I$ but often with
great deal fewer generators. We indicate  how this may come
about, with the  notion of {\em minimal reduction}.
Let $(R,{\frak m})$ be a Noetherian local ring
and let $I$ be an ideal (or a homogeneous ideal of a graded ring).
The  {\em special fiber} of the Rees algebra $R[It]$ is the ring
\[\F(I)=R[It]\otimes_R R/{\frak m}.\]
Its Krull dimension is called the {\em analytic spread} of $I$, and is
denoted $\ell(I)$.

If $R/{\frak m}$ is an infinite field, minimal
reductions of $I$ arise from the standard Noether normalizations of
the graded algebra $\F(I)$. The number of minimal generators of such
reductions is $\ell(I)$.
Let
\[
A= k[z_1, \ldots, z_{\ell}] \hookrightarrow \F(I),
\]
where $\ell = \ell(I) $, be a Noether normalization with the
$z_j\mbox{'}s$ chosen in degree $1$.
Let further $b_1, \ldots, b_s$ be a minimal set of homogeneous module
generators of $\F(I)$ over the algebra $A$
\[ \F(I) = \sum_{1\leq q \leq s}Ab_q. \]
If $J = (y_1, \ldots, y_{\ell})$, where $y_i$ is a lift in $R$ of
$z_i$, it is easy to see that $J$ is a reduction of $I$ and
$r_J(I)=\sup\{\deg b_q\}$.
In the case that the algebra $\F(I)$ is Cohen--Macaulay, $\F(I)$ is a
free module over $A$ so that $r_J(I)$ can be read off its
Hilbert--Poincar\'e series.

We shall now outline our results. In Section $2$, we relate the
exponent in the Dedekind--Mertens' formula directly to the
$a$-invariant of the Segre product of two rings of polynomials
(Theorem~\ref{sharpind}). The application of the formula to Noether
normalization is also pointed out in \cite{Edwards}. After remarks on
Gaussian ideals defined through algebras which are not polynomial
rings in Section $3$, we give in Section $4$ the primary decomposition
of the Gaussian ideal defined by two generic polynomials. The
components have the pleasing property that they all are Gorenstein
ideals (Theorem~\ref{decompofGfg}). In the final Section, we study the
normality of algebras associated to graphs; that includes the toric
algebras connected directly to (\ref{contentformula}). There are some
natural Noether normalizations for some of these extensions but not
the most general ones.

\section{Graphs and determinantal ideals}
If $G$ is a graph with vertices labelled by $x_0, \ldots, x_m$, its
monomial subring $k[G]$ is the subring of $k[x_0, \ldots, x_m]$
generated by all monomials $x_ix_j$ where $(x_i,x_j)$ is an edge of
$G$.
In parallel, there exists another algebra attached to $G$, defined by
the ideal of $k[x_0, \ldots, x_m]$ generated by those monomials (see
\cite{Vi}).
In general, it is difficult to find Noether normalizations of any of
these two families of algebras.

The following `explains' (\ref{contentformula2}) at the same time that
solves the question of Noether normalizations\footnote{After a first
draft of this note, we have found that \cite[Part 0]{Edwards} already
points out this Noether normalization. In addition, it has a delighful
historical account of (\ref{contentformula}). Our contribution on this
point is to explain the meaning of the exponent.}
for maximal bipartite graphs. It would be nice to find explicit
normalizations for other classes of graphs.

\begin{thm}\label{sharpind}
Let $X= \{x_0, \ldots, x_m\}$ and
$Y= \{y_0, \ldots, y_n\}$ be distinct sets of indeterminates and let
\[
f=\sum_{i=0}^m x_it^i \quad \mbox{and} \quad g=\sum_{j=0}^n y_jt^j
\]
be the corresponding generic polynomials over a field $k$. Set $R=k[X,
Y]$, $I = c(f)c(g)$, and $J= c(f g)$ and suppose $m\leq n$. Then
\begin{itemize}
\item[$($\mbox{a}$)$]
$J$ is a minimal reduction of $I$, $\ell(I)= m+n+1$, and $r_J(I) = m$.
\item[$($\mbox{b}$)$]
The polynomials
\[h_q= \sum_{i+j=q}x_iy_j\]
are algebraically independent and $k[h_q\mbox{\rm '}s]$ is a Noether
normalization of $k[x_iy_j\mbox{\rm '}s]$.
\end{itemize}
In particular, the factor $c(f)^m$ in the content formula~{\rm
(\ref{contentformula})} is sharp.
\end{thm}
\begin{pf}
We note that the ideal $I = (x_iy_j\mbox{\rm '}s)$ is the edge ideal
associated to the graph $G$ which is the join of two discrete graphs,
one with $m+1$ vertices and another with $n+1$ vertices; $G$ is,
therefore, bipartite.

Since $J$ is already a reduction of $I$ by (\ref{contentformula2}), we
may assume that $k$ is an infinite field.
On the other hand, as $I$ is generated by homogeneous polynomials of
the same degree, $\F(I) \simeq k[x_iy_j\mbox{\rm '}s]=k[G]$ (see
\cite{graphs}).
Let $Q_{ij}, 0\leq i\leq m,\ 0\leq j\leq n$ be distinct
indeterminates and map
\[ \psi: k[Q_{ij}\mbox{\rm '}s] \longrightarrow  k[x_iy_j\mbox{\rm '}s],\quad
\psi(Q_{ij})=x_iy_j.\]
We claim that the kernel of $\psi$ is generated by the $2 \times 2$ minors
of a generic $(m+1) \times (n+1)$ matrix.
Indeed let $\Q=(Q_{ij})$. It is clear that the ideal $I_2(\Q)$,
generated by the $2\times 2$ minors of $\Q$, is contained in ${\frak
Q}={\rm ker}(\psi)$. On the other hand, since the graph is bipartite,
$\dim(k[G])=m+n+1$ (see \cite{graphs}) and, therefore,
\[\mbox{\rm height}({\frak Q})= (m+1)(n+1)-(m+n+1)=m n=\mbox{\rm
height}(I_2(\Q)),\]
the latter by the classical formula for determinantal ideals
(see \cite[Theorem 2.5]{BV}).
Since they are both  prime ideals, we have $I_2(\Q)={\frak Q}$.

To complete the proof we note that the $a$-invariant of
$k[Q_{ij}\mbox{\rm '}s]/I_2(\Q)$ is $-n-1$ according to \cite{BHa}, and
therefore the reduction number of $\F(I)$ is $(m+n+1)-n-1=m$. \qed
\end{pf}

\begin{rem} {\rm
Another approach to the computation of the $a$-invariant is through
the theory of Segre products, and then appealing directly to \cite{GWI}.
The Cohen--Macaulayness and Gorensteiness of algebras that include
$k[x_iy_j\mbox{\rm '}s]$, has been dealt with in great detail already
in \cite{BST}. }
\end{rem}

\section{Generalized contents}
Let $R$ be a commutative ring and let $A$ be an $R$-algebra which is
free as an $R$-module. Let $\{e_{i}\mbox{\rm '}s\}$ be an $R$-basis with
attached structure constants $c_{ijk}$. Given an element $f\in A$,
define now $c(f)$ as the $R$-ideal generated by the coefficients of
the expression of $f$ as a linear combination of the $e_i\mbox{\rm '}s$.
This ideal is independent of the choice of basis being the usual
order ideal of an element of a free module.

We would like to know which condition on the $c_{ijk}\mbox{\rm '}s$
implies that $c(f g)$ is a reduction of $c(f) c(g)$.
Here is one instance

\begin{prop}
Let $A$ be an algebra which is a free module over the integral domain
$R$. Suppose $A$ has an $R$-basis indexed by a well-ordered monoid. If
for each pair of indices $i,j$ \ $\sum_kRc_{ijk}=R$, then for any two
elements $f,g\in A$, $c(f g)$ is a reduction of $c(f) c(g)$.
\end{prop}
\begin{pf}
We may replace $R$ by one of its valuation overrings $V$ (see
\cite[p. 350]{ZS}). It will then suffice to show that $c(f) c(g)V=c(f
g)V$.
\end{pf}

The assertion will follow from

\begin{lem}[{\bf Gauss Lemma}]
Let $A$ be an algebra as above and let $f$ and $g$ be two unimodular
elements $($i.e., $c(f)=c(g)=R)$. Then $f g$ is unimodular.
\end{lem}

\begin{rem}{\rm The condition on the well-ordering of the index
set of the basis is too restrictive, although it can be used for
bases change (for instance, even in the case of $R[t]$ one could use
other bases than   $\{t^n,\ n\geq 0\}$, with a
compatible ordering).
More precisely, once Gauss Lemma holds for a basis it will hold for
any other bases: all that requires is that for each prime $\frak p$ of
$R$ the fiber $A\otimes_R k({\frak p})$ is an integral domain.
}
\end{rem}

\section{Primary decomposition}
The generic form of the ideal $c(fg)$ has an interesting primary
decomposition.

\begin{thm}\label{decompofGfg}
Let $R$ be a Noetherian integral domain and let
\[ f = x_0 + x_1t + \cdots + x_mt^m \quad \mbox{and} \quad g = y_0 +
y_1t + \cdots + y_nt^n, \]
be generic polynomials of degrees $m$ and $n$ over $R$. The Gaussian ideal
$G(f, g)= c(fg)$ has a primary decomposition
\begin{equation}
c(fg) = c(f)\cap c(g) \cap [c(fg) + c(f)^{n+1} + c(g)^{m+1}].
\end{equation}
Furthermore, if $R$ is a Gorenstein ring then
\begin{equation}
L(f, g)=c(fg) + c(f)^{n+1} + c(g)^{m+1}
\end{equation}
is a Gorenstein ideal.
\end{thm}
\begin{pf}
The primary decomposition is easy to verify
\begin{eqnarray*}
c(f)\cap c(g) \cap [c(fg)+ c(f)^{n+1} + c(g)^{m+1}]   & = &        \\
c(f)\cap [c(g) \cap [c(fg)+ c(f)^{n+1} + c(g)^{m+1}]] & = &        \\
c(f) \cap  [c(fg)+ c(g)\cap c(f)^{n+1} + c(g)^{m+1}]  & = &        \\
c(f) \cap [c(fg) + c(g)^{m+1}]                        & = &        \\
c(fg) + c(f)\cap c(g)^{m+1}                           & = & c(fg),
\end{eqnarray*}
where in the third and fourth equalities we used
(\ref{contentformula}).

To prove that $L(f, g) = c(fg) + c(f)^{n+1} + c(g)^{m+1}$ is
Gorenstein, we show that it is a proper specialization of the
Gorenstein ideal described in \cite[Example 3.4]{HU}.

The building blocks of this ideal are a sequence ${\bf X} = (X_1,
\ldots, X_r)$ and a $r \times s$ matrix $\varphi$, $s\geq r$. For the
generic sequence and matrix define
\[ J = ({\bf X}\cdot \varphi) + ({\bf X})^{s-r+1} + I_{r}(\varphi), \]
where $({\bf X}\cdot \varphi) $ denotes the ideal generated by the
entries of the product of the sequence by the matrix, and
$I_{r}(\varphi)$ is the ideal generated by the minors of order $r$ of
the matrix $\varphi$. In \cite{HU} it is shown that $J$ is a
Gorenstein ideal of codimension $s+1$.

In our case,
\begin{eqnarray*}
{\bf X} & = & (x_0, x_1, \ldots, x_m)
\end{eqnarray*}
and $\varphi$ is the $(m+1) \times (m+n+1)$ matrix
\[  \varphi = \left[ \begin{array}{cccccccc}
y_0 & y_1 & y_2 & \cdots & y_n & 0 & \cdots  & 0 \\
0   & y_0 & y_1 & \cdots & y_{n-1} & y_n & \cdots & 0 \\
\vdots & \vdots & \vdots & \vdots & \vdots & \vdots & \vdots & \vdots \\
0 & 0 & \cdots & y_0 & y_1 & y_2 & \cdots & y_n
\end{array} \right]. \]

Now we note that
\begin{eqnarray*}
({\bf X}\cdot \varphi)  & = & c(fg), \\
({\bf X})^{s-r+1}       & = & c(f)^{n+1}, \\
I_r(\varphi)            & = & c(g)^{m+1},
\end{eqnarray*}
and the codimension of $L(f,g)$ is $s+1 = m+n+2$. This means that
$L(f,g)$ is a proper specialization of $J$ and is therefore Gorenstein
as well. \qed
\end{pf}

It is natural to define the Gaussian ideal associated to any finite
set of polynomials. In the generic case, these ideals share similar
properties to $G(f,g)$. Let us consider the case of $3$ polynomials,
where an open question arises.

\begin{thm}\label{prim-dec2}
Let $X=\{x_0, \ldots, x_m\}$, $Y=\{y_0, \ldots,y_n\}$, and $Z=\{z_0,
\ldots, z_p\}$ be $3$ sets of indeterminates. Defining the polynomials
\[ f=\sum_{i=0}^mx_it^i, \quad g=\sum_{j=0}^ny_jt^j, \quad \mbox{and}
\quad h=\sum_{k=0}^pz_kt^k, \]
one has that the primary decomposition of $c(f g h)$ is given by
\begin{equation}
c(fgh)= c(f) \cap c(g) \cap c(h) \cap L(f, g) \cap L(f, h) \cap L(g,
h) \cap L(f, g, h),
\end{equation}
where
\begin{eqnarray*}
L(f, g)       & = & c(fg) + c(f)^{n+1} + c(g)^{m+1}, \\
L(f, h)       & = & c(fh) + c(f)^{p+1} + c(h)^{m+1}, \\
L(g, h)       & = & c(gh) + c(g)^{p+1} + c(h)^{n+1}, \\
L(f, g, h)    & = & c(fgh) + c(fg)^{p+1} + c(fh)^{n+1} + c(gh)^{m+1} + \\
              & \ & \mbox{\hspace{1.5cm}} c(f)^{n+p+1} + c(g)^{m+p+1} +
                    c(h)^{m+n+1}.
\end{eqnarray*}
\end{thm}
\begin{pf}
The proof follows from a repeated use of Theorem~\ref{decompofGfg} and
the Dedekind--Mertens formula. Indeed, one easily verifies that
\begin{eqnarray*}
& \ & \mbox{\hspace{-1cm}}c(f) \cap c(g) \cap c(h) \cap L(f, g) \cap
      L(f, h) \cap L(g, h) \cap L(f, g, h) = \\
& \ & \ \\
& = & c(fg) \cap c(fh) \cap c(gh) \cap L(f, g, h) \\
& = & c(fg) \cap c(fh) \cap \left[ \right. c(fgh) + c(gh) \cap
      c(fg)^{p+1} + c(gh) \cap c(fh)^{n+1} + \\
& \ & \mbox{\hspace{1cm}}
      c(gh)^{m+1} + c(gh) \cap c(f)^{n+p+1} + c(gh) \cap
      c(g)^{m+p+1} + \\
& \ & \mbox{\hspace{2cm}}
      c(gh) \cap c(h)^{m+n+1} \left. \right] \\
& = & c(fg) \cap c(fh) \cap [ c(fgh) + c(gh) \cap c(fg)^{p+1} + c(gh)
      \cap c(fh)^{n+1} + \\
& \ & \mbox{\hspace{1cm}}
      c(gh)^{m+1} + c(fgh)c(f)^{n+p} + c(gh) \cap c(g)^{m+p+1} + \\
& \ & \mbox{\hspace{2cm}}
      c(gh) \cap c(h)^{m+n+1} ] \\
& = & c(fg) \cap [ c(fgh) + c(fh) \cap c(gh) \cap c(fg)^{p+1} + c(gh)
      \cap c(fh)^{n+1} + \\
& \ & \mbox{\hspace{1cm}}
      c(fh) \cap c(gh)^{m+1} + c(gh) \cap c(fh) \cap
      c(g)^{m+p+1} + \\
& \ & \mbox{\hspace{2cm}}
      c(fh) \cap c(gh) \cap c(h)^{m+n+1} ] \\
& = & c(fg) \cap [ c(fgh) + c(fh) \cap c(gh) \cap c(fg)^{p+1} +
      c(gh) \cap c(fh)^{n+1} + \\
& \ & \mbox{\hspace{1cm}}
      c(fh) \cap c(gh)^{m+1} + c(gh) \cap
      c(fgh)c(g)^{m+p} + \\
& \ & \mbox{\hspace{2cm}}
      c(fh) \cap c(gh) \cap c(h)^{m+n+1} ] \\
& = & c(fgh) + c(fh) \cap c(gh) \cap c(fg)^{p+1} + c(fg) \cap c(gh) \cap
      c(fh)^{n+1} + \\
& \ & \mbox{\hspace{1cm}}
      c(fg) \cap c(fh) \cap c(gh)^{m+1} + c(fh)
      \cap c(gh) \cap c(fg) \cap c( h)^{m+n+1} \\
& = & c(fgh) + c(fh) \cap c(gh) \cap c(fg)^{p+1} + c(fg) \cap c(gh) \cap
      c(fh)^{n+1} + \\
& \ & \mbox{\hspace{1cm}}
      c(fg) \cap c(fh) \cap c(gh)^{m+1} + c(fh)
      \cap c(gh) \cap c(fgh)c(h)^{m+n} \\
& = & c(fgh) + c(h) \cap c(fg)^{p+1} + c(g)
      \cap c(fh)^{n+1} + c(f) \cap c(gh)^{m+1}  \\
& = & c(fgh) + c(fgh) \cap c(fg)^{p} + c(fgh)
      \cap c(fh)^{n} + c(fgh) \cap c(gh)^{m} \\
& = & c(fgh),
\end{eqnarray*}
as claimed. \qed
\end{pf}

\begin{rem} {\rm
Experiments show that the ideals $L(f,g,h)$ are Gorenstein. Perhaps
they can be obtained by specialization of sums of Huneke--Ulrich ideals. }
\end{rem}

\section{Multiproducts and joins}
In order to see a different explanation of (\ref{contentformula2}), we
extend it to the product of 3 (or more) polynomials, but use the theory
of Segre products as a tool.

Let $X=\{x_0, \ldots, x_m\}$, $Y=\{y_0,\ldots,y_n\}$, and $Z=\{z_0,\ldots,
z_p\}$ be 3 sets of indeterminates. Defining the polynomials
\[ f=\sum_{i=0}^m x_it^i, \quad g=\sum_{j=0}^n y_jt^j, \quad \mbox{and} \quad
h=\sum_{k=0}^p z_kt^k, \]
one has that $J=c(f g h)$ is a reduction of $I=c(f)c(g)c(h)$ by
Gauss Lemma. If $m\leq n\leq p$, a simple calculation will show that
$\ell(I) = m+n+p+1$ and $r_J(I)\leq m+n$. We now resolve this inequality.

\begin{prop}
Let $X=\{x_0, \ldots, x_m\}$, $Y=\{y_0, \ldots, y_n\}$, and
$Z=\{z_0, \ldots, z_p\}$ be sets of distinct indeterminates, let
$R=k[X,Y,Z]$ be a polynomial ring over a field $k$, and let
$I=(x_iy_jz_k\vert\, x_i\in X,\, y_j\in Y,\, z_k\in Z)$. Then $I$ is a
normal ideal of $R$.
\end{prop}
\begin{pf}
We will show that $I^q$ is complete for all
$q\geq 1$. Let $I_a^q$ be the integral closure of $I^q$ and let
$f\in I_a^q$ be a monomial. We write
\[ f=x_{i_1}^{a_1}\cdots x_{i_r}^{a_r}y_{j_1}^{b_1} \cdots
y_{j_s}^{b_s} z_{k_1}^{c_1}\cdots z_{k_t}^{c_t}. \]
Since $f^w\in I^{sq}$ for some $w>0$ we can write
\[f^w=x_{q_1}^{d_1}\cdots
x_{q_\lambda}^{d_\lambda}M,
\]
where $M$ is a monomial whose support is contained in $Y\cup Z$. We
obtain $w\sum_{i=1}^ra_i=\sum_{i=1}^{\lambda}d_i\geq wq$, which
implies $\sum_{i=1}^ra_i\geq q$, and a similar argument shows
$\sum_{i=1}^sb_i\geq q$ and $\sum_{i=1}^tc_i\geq q$. Therefore $f\in
I^q$. \qed
\end{pf}

Note that by Hochster's theorem (see \cite{Hochster}), the
algebra $R[It]$ is Cohen--Macaulay. Furthermore, since
$\F(I)= k[x_iy_jz_k\vert\, x_i\in X,\, y_j\in Y,\, z_k\in Z]$ is a
direct summand of $R[It]$, it is also normal and therefore
Cohen--Macaulay by \cite{Hochster}. We may thus more easily compute
the reduction number of $\F(I)$.

\begin{thm}
The reduction number of the ideal $I$ above is $m+n$.
\end{thm}
\begin{pf}
Since $\F(I)$ is Cohen--Macaulay, its reduction number can also
be obtained from the degrees of the generators of its canonical
module. But $\F(I)$ is a Segre product of standard Cohen--Macaulay
algebras and the canonical module is given by an explicit formula
from the canonical modules of the factors (see \cite[Theorem
4.3.1]{GWI}). Entering the data we get $r_J(I) = m+n$. \qed
\end{pf}

\begin{rem}
{\rm Semigroup rings attached to more general bipartite graphs are
obtained by deleting some of the generators in
$k[x_iy_j\mbox{\rm '}s]$. These rings are still normal but we do not know
what their  reductions are like. }
\end{rem}

\subsection*{The join of two normal ideals}
\begin{defn} {\rm
Let $I$ and $J$ be two monomial ideals of the polynomial rings
$k[x_0,\ldots,x_m]$ and $k[y_0,\ldots, y_n]$ respectively. The {\em
join} of $I$ and $J$ is:
\[ I\ast J=I+J+K;\ \mbox{where } K=(x_iy_j\vert\, 0\leq i \leq m\mbox{
and } 0\leq j\leq n). \] }
\end{defn}

\begin{thm}\label{join}
Let $R=k[x_0,\ldots,x_m]$ and $S=k[y_0,\ldots,y_n]$ be polynomial rings
over a field $k$, and let $I$, $J$ be two ideals of $R$ and $S$
respectively. If $I$ and $J$ are normal ideals generated by
square-free monomials of the same degree $t\geq 2$ then their join
$I\ast J$ is normal.
\end{thm}
\begin{pf}
Set $X=\{x_0,\ldots,x_m\}$, $Y=\{y_0,\ldots,y_n\}$ and
$L=I+J+K$, where $K=(X)(Y)$. By induction on $p$ we will show that
$L_a^p=L^p$ for all $p\geq 1$, where $L_a^p$ denotes the integral
closure of $L^p$. If $p=1$ then $L$ is a radical ideal
(see \cite[Prop.~1]{EH}), hence $L$ is integrally closed. Assume
$L_a^i=L^i$ for  $i<p$ and $p\geq 2$. Using the results of
\cite{Kempf} we have
\[
L_a^p=(\{z\vert\, z\mbox{ is a monomial in }k[X,Y]\mbox{  and }z^q\in
L^{qp}\mbox{ for  some }q\geq 1\}).
\]
Let $z$ be a monomial in $L_a^p$, then $z^q\in L^{qp}$, $q>0$. Let us
show $z\in L^p$. Since $L_a^p\subseteq
L_a^{p-1}=L^{p-1}$ we can write
\[
z=Mh_1\cdots h_sg_1\cdots g_rf_{1}\cdots
f_{p-r-s-1},
\]
where $M$ is a monomial, the $h_i$'s are monomials of
degree two in $K$, the $g_i$'s and $f_i$'s are degree $t$ monomials
in $J$ and $I$ respectively. Likewise we can write
\[
z^q=Nh_1'\cdots h_{s_1}'g_1'\cdots g_{r_1}'f_{1}'\cdots
f_{qp-r_1-s_1}',
\]
where $N$ is a monomial, $\deg(h_i')=2$ and $h_i'$ is a monomial in
$K$ for all $i$, $g_i'$ and $f_j'$ are degree $t$ monomials
in $J$ and $I$ respectively for all $i,j$. From the last two
equalities we have
\begin{eqnarray}\label{from1-2}
z^q &=&M^q(h_1\cdots h_s)^q(g_1\cdots g_r)^q(f_{1}\cdots
f_{p-r-s-1})^q\\
 &=& Nh_1'\cdots h_{s_1}'g_1'\cdots g_{r_1}'f_{1}'\cdots
f_{qp-r_1-s_1}'.\nonumber
\end{eqnarray}
From (\ref{from1-2}) one readily derives the inequality
\begin{equation}\label{from1-2-3}
(qs-s_1)(t-2)+qt\leq q\deg(M).
\end{equation}
We may assume $M=x^\alpha$ or $M=y^\alpha$, otherwise $z\in L^p$. By
symmetry we may assume  $M=x^\alpha$, where $\alpha\geq 0$.

({\em a}) If $t\geq 3$, then $r=0$ or $p-r-s-1=0$, otherwise $g_1f_1\in
L^3$ and hence $z\in L^p$. First we treat the case $r=0$. By taking
degrees in (\ref{from1-2}) w.r.t. the variables $y_0,\ldots,y_n$ one
has $s_1+tr_1\leq qs$, which together with (\ref{from1-2-3}) yields
$\deg M\geq t$. Let $z_1$ be the result of evaluating $z$ at $y_i=1$.
From (\ref{from1-2}) we derive
$z_1^q\in I^{qp-s_1-r_1}\subseteq I^{q(p-s)}$, and $z_1\in
I_a^{p-s}=I^{p-s}$. We may write $z=y^\theta z_1$ and $z_1=x^\beta w$,
where $\deg(y^\theta)=s$ and $w$ is a monomial in $I^{p-s}$ of degree
$t(p-s)$. Since $\deg(z_1)=\deg(M)+s+t(p-s-1)$ we obtain
$\deg(z_1)\geq s+t(p-s)$, hence $\deg(x^\beta)\geq s$. Altogether
we derive $z=y^\theta x^\beta w\in L^p$. Next we consider the case
two $r=p-s-1\geq 1$, observe that $\deg(M)\leq 1$, otherwise
$z\in L^p$. Therefore either $z=h_1\cdots h_s g_1\cdots g_r$, or we
may rewrite $z=y^\beta h_1\cdots
h_sh_{s+1}g_1\cdots g_{r-1}$, where $\deg(h_{s+1})=2$ and $h_{s+1}\in
K$, interchanging the $x_i$ and $y_i$ variables we may apply the
arguments above to conclude $z\in L^p$.

({\em b}) Assume $t=2$. Using $z^q\in L^{qp}$ one rapidly obtains
$\deg(M)\geq 2$, hence we may assume $r=0$ (otherwise
$z\in L^p$) and the arguments of case ({\em a}) can be applied to conclude
$z\in L^p$. \qed
\end{pf}

The following Corollary generalizes the normality assertion of
\cite[Theorem~4.8(v)]{aron}.

\begin{cor}
Let $X=\{x_0,\ldots,x_m\}$ and $\{y_0,\ldots,y_n\}$ be two disjoint
sets of indeterminates over a field $k$. Let $I$ be a normal ideal of
$k[X]$
generated by
square free monomials of degree $t$ and let $L=I+K$, where
$K=(X)(Y)$. Then $L$ is a normal ideal.
\end{cor}
\begin{pf}
Proceed as in the proof of Theorem~\ref{join} and notice that
in this case $r=p-s-1$. \qed
\end{pf}

\begin{ack}
A. Corso gratefully acknowledges partial support from the {\it Consiglio
Nazionale delle Ricerche} under CNR grant 203.01.63. W. V. Vasconcelos
was in part supported by a NSF grant. Finally, R. H. Villarreal expresses his
gratitude to COFAA-IPN, CONACyT, and SNI.
\end{ack}

\end{document}